\documentclass[11pt]{amsproc}%
\usepackage{amsfonts}
\usepackage{amsmath}
\usepackage{amssymb}
\usepackage{graphicx}%
\renewcommand {\theequation} {\@arabic\c@equation}
 \theoremstyle{plain}

\newtheorem{corollary}{Corollary}

\newtheorem{definition}{Definition}

\newtheorem{lemma}{Lemma}

\newtheorem{proposition}{Proposition}
\newtheorem{remark}{Remark}

\newtheorem{theorem}{Theorem}
\numberwithin{equation}{section}
\begin{document}

\begin{center}

\end{center}

\centerline{\Large {\bf Free Orbit-Dimension of Finite von Neumann
Algebras}}

$\bigskip$

\centerline{\large Don Hadwin \qquad and \qquad Junhao Shen}

\bigskip

\centerline{Mathematics Department, University of New Hampshire,
Durham, NH, 03824}

\bigskip

\centerline{email: don@math.unh.edu \qquad and \qquad
jog2@cisunix.unh.edu}

$\bigskip$

\noindent\textbf{Abstract: } We introduce a new free entropy invariant, which
  yields significant improvements of most of the applications of free
entropy to finite von Neumann algebras, including those with Cartan
subalgebras, simple masas, property $T,$ property $\Gamma,$ nonprime factors,
and thin factors.

\vspace{0.3cm}

\section{Introduction}

The theory of free probability and free entropy was introduced by
Voiculescu in 1980's.  In his  papers \cite{V2} \cite{V3},
Voiculescu introduced the concept of free entropy dimension and used
it to provide the first example of II$_{1}$ factor that does not
have Cartan subalgebras, which solves a long-standing open problem.
 Later Ge in \cite{Ge2} showed that the free group factors   are
  not prime, i.e., are not a tensor product of two infinite-dimensional von Neumann algebras. This   also answers a very old open
  question.
 In \cite{GS2}, Ge and the second author computed free entropy
dimension for a large class of finite von Neumann algebras including
some II$_{1}$ factors with property $T$.

Here we introduce a new invariant, the upper free orbit-dimension of
a finite von Neumann algebra, which is closely related to
Voiculescu's free entropy dimension. Suppose that $\mathcal M$ is a
von Neumann algebra with a tracial state $\tau$. Roughly speaking,
if $x_{1},\ldots,x_{n}$ generates $\mathcal{M},$ Voiculescu's free
entropy dimension $\delta_{0}\left( x_{1},\ldots ,x_{n}\right)  $ is
obtained by considering the covering numbers of certain sets by
$\omega$-balls, and letting $\omega$ approach $0$. The upper free
orbit-dimension $\mathfrak{K}_{2}(x_{1},\ldots,x_{n})$ is obtained
by considering the covering numbers of the same sets by
$\omega$-neighborhoods of unitary orbits (see the definitions in
section 2), and taking the supremum over $\omega,$ $0<\omega<1$. It
is easily shown that
\[
\delta_{0}\left(  x_{1},\ldots,x_{n}\right)  \leq1+\mathfrak{K}_{2}%
(x_{1},\ldots,x_{n})
\]
always \ holds. \ Most \ of \ the \ important \ applications \
involving \ $\delta_{0}$ \ involve \  showing \\ $\delta_{0}\left(
x_{1},\ldots,x_{n}\right)  \leq1,$ while we see that it is much
easier to show $\mathfrak{K}_{2}(x_{1},\ldots,x_{n})=0$.

The upper free orbit-dimension has many useful properties, mostly in the case
when $\mathfrak{K}_{2}(x_{1},\ldots,x_{n})=0.$ The key property is that if
$\mathfrak{K}_{2}(y_{1},\ldots,y_{p})=0$ for some generating set for
$\mathcal{M},$ then $\mathfrak{K}_{2}(x_{1},\ldots,x_{n})=0$ for every
generating set. This fact allows us to show that the class of finite von
Neumann algebras $\mathcal{M}$ with $\mathfrak{K}_{2}(\mathcal{M})=0$ is
closed under certain operations that enlarge the algebra:

\begin{enumerate}
\item If $\mathfrak{K}_{2}(\mathcal{N}_{1})=\mathfrak{K}_{2}(\mathcal{N}%
_{2})=0$ and $\mathcal{N}_{1}\cap\mathcal{N}_{2}$ is diffuse, then
$\mathfrak{K}_{2}(\left(  \mathcal{N}_{1}\cup\mathcal{N}_{2}\right)
^{\prime\prime})=0.$

\item If $\ \mathcal{M}=\{\mathcal{N},u\}^{\prime\prime}$ where $\mathcal{N}$
is a von Neumann subalgebra of $\ \mathcal{M}$ with $\mathfrak{K}%
_{2}(\mathcal{N})=0$ and $u$ is a unitary element in $\ \mathcal{M}$
satisfying, for a sequence $\left\{  v_{n}\right\}  $ of Haar unitary elements
in $\mathcal{N}$, dist$_{\left\Vert {}\right\Vert _{2}}\left(  uv_{n}u^{\ast
},\mathcal{N}\right)  \rightarrow0$, then $\mathfrak{K}_{2}(\ \mathcal{M})=0$.
\item If $\{\mathcal{N}_{i}\}_{i=1}^{\infty}$ is an ascending sequence of von
Neumann subalgebras of $\mathcal{M}$ such that $\mathfrak{K}_{2}%
(\mathcal{N}_{i})=0$ for all $i\geq1$ and $\mathcal{M}=\overline{\cup
_{i}\mathcal{N}_{i}}^{SOT}$, then $\mathfrak{K}_{2}(\mathcal{M})=0$. \bigskip
\end{enumerate}
  Using these closure operations as building blocks, and the
easily-proved fact that $\mathfrak{K}_{2}(\mathcal{M})=0$ whenever
$\mathcal{M}$ is hyperfinite, we can can show that
$\mathfrak{K}_{2}(\mathcal{M})=0$ for a large class of von Neumann
algebras. As a corollary we  recapture most of the old results. In
particular, we
  extend results in \cite{V3},
\cite{Ge2}, \cite{GS2}, \cite{GS1},  \cite{GePopa}, \cite{H},
\cite{V4},  \cite{Ge1}, \cite{Dyk}.

\section{Definitions}

Let $\mathcal{M}_{k}(\mathbb{C})$ be the $k\times k$ full matrix
algebra with entries in $\mathbb{C}$, and $\tau_{k}$ be the
normalized trace on $\mathcal{M}_{k}(\mathbb{C})$, i.e.,
$\tau_{k}=\frac{1}{k}Tr$, where $Tr$ is the usual trace on
$\mathcal{M}_{k}(\mathbb{C})$. Let $\mathcal{U}(k)$ denote the group
of all unitary matrices in $\mathcal{M}_{k}(\mathbb{C})$. Let
$\mathcal{M}_{k}(\mathbb{C})^{n}$ denote the direct sum of $n$
copies of $\mathcal{M}_{k}(\mathbb{C})$. Let $\Vert\cdot\Vert_{2}$
denote the trace norm induced by $\tau_{k}$ on
$\mathcal{M}_{k}(\mathbb{C})^{n}$, i.e.,
\[
\Vert(A_{1},\ldots,A_{n})\Vert_{2}^{2}=\tau_{k}(A_{1}^{\ast}A_{1})+\ldots
+\tau_{k}(A_{n}^{\ast}A_{n})
\]
for all $(A_{1},\ldots,A_{n})$ in $\mathcal{M}_{k}(\mathbb{C})^{n}$.

For every $\omega>0$, we define the $\omega$-ball $Ball(B_{1},\ldots
,B_{n};\omega)$ centered at $(B_{1},\ldots,B_{n})$ in
$\mathcal{M}_{k}(\mathbb{C})^{n}$ to be the subset of
$\mathcal{M}_{k}(\mathbb{C})^{n}$   consisting of all
$(A_{1},\ldots,A_{n})$ in $\mathcal{M}_{k}(\mathbb{C})^{n}$ such
that
$\Vert(A_{1},\ldots,A_{n})-(B_{1},\ldots,B_{n})\Vert_{2}<\omega.$

For every $\omega>0$, we define the $\omega$-orbit-ball $\mathcal{U}%
(B_{1},\ldots,B_{n};\omega)$ centered at $(B_{1},\ldots,B_{n})$ in $\mathcal{M}%
_{k}(\mathbb{C})^{n}$ to be the subset of
$\mathcal{M}_{k}(\mathbb{C})^{n}$ consisting of all
$(A_{1},\ldots,A_{n})$ in $\mathcal{M}_{k}(\mathbb{C})^{n}$ such
that there exists some unitary matrix $W$ in $\mathcal{U}(k)$
satisfying
\[
\Vert(A_{1},\ldots,A_{n})-(WB_{1}W^{\ast},\ldots,WB_{n}W^{\ast})\Vert
_{2}<\omega.
\]

Let $\mathcal{M}$ be a von Neumann algebra with a tracial state
$\tau$, and $x_{1},\ldots,x_{n}$ be elements in $\mathcal{M}$. We
now define our new invariants. For any positive $R$ and $\epsilon$,
and any $m,k$ in $\mathbb{N}$, let
$\Gamma_{R}(x_{1},\ldots,x_{n};m,k,\epsilon)$ be the subset of $\mathcal{M}%
_{k}(\mathbb{C})^{n}$ consisting of all $(A_{1},\ldots,A_{n})$ in
$\mathcal{M}_{k}(\mathbb{C})^{n}$ such that $\Vert A_{j}\Vert\leq R$, $1\leq
j\leq n$, and
\[
|\tau_{k}(A_{i_{1}}^{\eta_{1}}\cdots A_{i_{q}}^{\eta_{q}})-\tau(x_{i_{1}%
}^{\eta_{1}}\cdots x_{i_{q}}^{\eta_{q}})|<\epsilon,
\]
for all $1\leq i_{1},\ldots,i_{q}\leq n$, all $\eta_{1},\ldots,\eta_{q}$ in
$\{1,\ast\}$, and all $q$ with $1\leq q\leq m$.

For $\omega>0$, we define the $\omega$-orbit covering number $\nu(\Gamma
_{R}(x_{1},\ldots,x_{n};m,k,\epsilon),\omega)$ to be the minimal number of
$\omega$-orbit-balls that cover $\Gamma_{R}(x_{1},\ldots,x_{n};m,k,\epsilon)$
with the centers of these $\omega$-orbit-balls in $\Gamma_{R}(x_{1}%
,\ldots,x_{n};m,k,\epsilon)$. Now we define, successively,
\[
\begin{aligned}
\frak K (x_1,,\ldots,
x_n;\omega,R) &=   \inf_{m\in \Bbb N,
\epsilon>0}\limsup_{k\rightarrow \infty} \frac
{\log(\nu(\Gamma_R(x_1,\ldots,
x_n;m,k,\epsilon),\omega))}{-k^2\log\omega} \\
\frak K (x_1,,\ldots,
x_n;\omega) &= \sup_{R>0}  \frak K (x_1,,\ldots,
x_n;\omega,R)\\
\frak K_1(x_1,,\ldots, x_n ) &= \limsup_{\omega\rightarrow 0}\frak K
(x_1,,\ldots,
x_n;\omega)\\
\frak K_2(x_1,,\ldots, x_n ) &= \sup_{0<\omega<1}\frak K
(x_1,,\ldots, x_n;\omega),
\end{aligned}
\]
where  $\mathfrak{K}_{1}(x_{1},,\ldots,x_{n})$   is called  the {\em
free
orbit-dimension} of \ $x_{1},\ldots,x_{n}$   and  $\mathfrak{K}%
_{2}(x_{1},,\ldots,x_{n})$ is called the  {\em upper free
orbit-dimension} of $x_{1},\ldots x_{n} $.

In the   spirit as in Voiculescu's definition of free entropy
dimension, we shall also define free orbit-dimension and upper free
orbit-dimension of $x_{1},\ldots,x_{n}$ in the presence of
$y_{1},\ldots,y_{p}$ for all \
$x_{1},\ldots,x_{n},y_{1},\ldots,y_{p}$ \ in \ the \ von \ Neumann \
algebra
\ $\mathcal{M}$ \ as \ follows. Let $\Gamma_{R}(x_{1},\ldots,x_{n}:y_{1}%
,\ldots,y_{p};m,k,\epsilon)$ be the image of the projection of
$\Gamma _{R}(x_{1},\ldots,x_{n},y_{1},\ldots,y_{p};m,k,\epsilon)$
onto the first $n$ components, i.e.,
\[
(A_{1},\ldots,A_{n})\in\Gamma_{R}(x_{1},\ldots,x_{n}:y_{1},\ldots
,y_{p};m,k,\epsilon)
\]
if there are elements $B_{1},\ldots,B_{p}$ in $\mathcal{M}_{k}(\mathbb{C})$
such that
\[
(A_{1},\ldots,A_{n},B_{1},\ldots,B_{p})\in\Gamma_{R}(x_{1},\ldots,x_{n}%
,y_{1},\ldots,y_{p};m,k,\epsilon).
\]
Then we define, successively,
\[
\begin{aligned}
\frak K (x_1,&\ldots,
x_n:y_1,\ldots,y_p;\omega,R)\\ & \qquad \qquad =   \inf_{m\in \Bbb
N, \epsilon>0}\limsup_{k\rightarrow \infty} \frac
{\log(\nu(\Gamma_R(x_1,\ldots,
x_n:y_1,\ldots,y_p;m,k,\epsilon),\omega))}{-k^2\log\omega} \\
\frak K (x_1,&\ldots, x_n:y_1,\ldots,y_p;\omega ) = \sup_{R>0} \frak
K (x_1, \ldots,
x_n:y_1,\ldots,y_p;\omega,R)\\
\frak K_1(x_1,&\ldots, x_n:y_1,\ldots,y_p ) =
\limsup_{\omega\rightarrow 0}\frak K (x_1,\ldots,
x_n:y_1,\ldots,y_p;\omega)\\
\frak K_2(x_1,&\ldots, x_n :y_1,\ldots,y_p)  =
\sup_{0<\omega<1}\frak K (x_1,\ldots, x_n:y_1,\ldots,y_p;\omega).
\end{aligned}
\]

\begin{definition}
Suppose $\mathcal{M}$ is a finitely generated von Neumann algebra with a
tracial state $\tau$. Then the \emph{free orbit-dimension} $\mathfrak{K}%
_{1}(\mathcal{M})$ of $\mathcal{M}$ is defined by
\[
\mathfrak{K}_{1}(\mathcal{M})=\sup\{\mathfrak{K}_{1}(x_{1},\ldots
,x_{n})\ |\ \text{$x_{1},\ldots,x_{n}$ generate }\mathcal{M}\text{ as a von
Neumann algebra}\},
\]
and the \emph{upper free orbit-dimension} $\mathfrak{K}_{2}(\mathcal{M})$ of
$\mathcal{M}$ is defined by
\[
\mathfrak{K}_{2}(\mathcal{M})=\sup\{\mathfrak{K}_{2}(x_{1},\ldots
,x_{n})\ |\ \text{$x_{1},\ldots,x_{n}$ generate }\mathcal{M}\text{ as a von
Neumann algebra}\},
\]

\end{definition}

Here, we quote a useful proposition from \cite{DH}
\begin{proposition}
Suppose $\mathcal M$ is a hyperfinite von Neumann algebra with a
tracial state $\tau$. Suppose that $x_1,\ldots,x_n$ is a family of
generators of $\mathcal M$. Then, for every $\omega>0$,
$R>\max_{1\le j\le n}\|x_j\|$, there are a positive integer $m$ and
a positive number $\epsilon$ such that the following hold: for all
$k\ge 1$, if $A_1,\ldots, A_n, B_1,\ldots, B_n$ in $\mathcal
M_k(\Bbb C)$ satisfying, (a) $ 0\le \|A_j\|, \|B_j\|\le R$ for all
$1\le j\le n$; (b)
$$
\begin{aligned}
 &|\tau_k(A_{i_1}^{\eta_1}\cdots A_{i_p}^{\eta_p})- \tau (x_{i_1}^{\eta_1}\cdots
 x_{i_p}^{\eta_p})|<\epsilon\\
 & |\tau_k(B_{i_1}^{\eta_1}\cdots B_{i_p}^{\eta_p})- \tau (x_{i_1}^{\eta_1}\cdots
 x_{i_p}^{\eta_p})|<\epsilon,
\end{aligned}
$$
 for all $1\le i_1,\ldots, i_p\le n,$ $ \{\eta_{j}\}_{j=1}^p\subset
 \{*, 1\}$ and $1\le p\le m$, then there exists a unitary matrix $U$
 in $\mathcal U(k)$ such that
 $$
\sum_{j=1}^n \|U^*A_jU-B_j\|_2 <\omega.
 $$

\end{proposition}

\section{Key Properties of $\mathfrak{K}_{2}$}

\begin{lemma}
Let $x_{1},\ldots,x_{n}$ be self-adjoint elements in a von Neumann
algebra $\mathcal{M}$ with a tracial state $\tau$. Let
$\delta_{0}(x_{1},\ldots ,x_{n})$ be  Voiculescu's   free entropy
dimension. Then
\[
\delta_{0}(x_{1},\ldots,x_{n})\leq\mathfrak{K}_{1}(x_{1},\ldots,x_{n}%
)+1\leq\mathfrak{K}_{2}(x_{1},\ldots,x_{n})+1.
\]

\end{lemma}

\begin{proof}
The first inequality follows from Theorem 14 in \cite{DH}, and the second
inequality is obvious.
\end{proof}


\begin{lemma}
Let $x_{1},\ldots,x_{n},y_{1},\ldots,y_{p}$ be elements in a von Neumann
algebra $\mathcal{M}$ with a tracial state $\tau$. If $y_{1},\ldots,y_{p}$ are
in the von Neumann subalgebra generated by $x_{1},\ldots,x_{n}$ in
$\mathcal{M}$, then, for every $0<\omega<1,$
\[
\mathfrak{K}(x_{1},\ldots,x_{n};\omega)=\mathfrak{K}(x_{1},\ldots,x_{n}%
:y_{1},\ldots,y_{p};\omega)\mathfrak{.}%
\]

\end{lemma}
\begin{proof}
\textbf{\ } It is a straightforward adaptation of the proof of Prop.
1.6 in \cite{V3} (see also Lemma 5 in \cite{DH}). Given
$R>\max_{1\leq j\leq
p}\Vert y_{j}\Vert $, $m\in \mathbb{N}$ and $\epsilon >0$, we can find $%
m_{1}\in \mathbb{N}$ and $\epsilon _{1}>0$ such that, for all $k\in \mathbb{N%
}$,$$
\begin{aligned}
\Gamma _{R}(x_{1},\ldots ,x_{n};m_{1},k,\epsilon _{1})&\subset    \
\Gamma_{R}(x_{1},\ldots ,x_{n}:y_{1},\ldots ,y_{p};m,k,\epsilon )\\
 & \subset \ \Gamma _{R}(x_{1},\ldots ,x_{n};m,k,\epsilon ).
\end{aligned}$$
 Hence
\[
\begin{aligned}
\nu(\Gamma_R(x_1,\ldots,x_n;m_1,k,\epsilon_1),\omega) &\le\nu(
\Gamma_R (x_1,\ldots,x_n: y_1,\ldots,y_p;m,k,\epsilon),\omega)\\
&  \le\nu(\Gamma_R(x_1,\ldots,x_n;m,k,\epsilon),\omega),
\end{aligned}
\]%
for all $0<\omega <1$. The rest follows from the definitions.
\end{proof}

\vspace{0.2cm}

\noindent The following key theorem shows that, in some cases, the
upper free orbit-dimension $\mathfrak{K_{2}}$ is a von Neumann
algebra invariant, i.e., it is independent of the choice of
generators.

\begin{theorem}
Suppose $\mathcal{M}$ is a von Neumann algebra with a tracial state
$\tau$ and is generated by a family of elements
$\{x_{1},\ldots,x_{n}\}$  as a von Neumann algebra. If
\[
\mathfrak{K}_{2}(x_{1},\ldots,x_{n})=0,
\]
then
\[
\mathfrak{K}_{2}(\mathcal{M})=0.
\]

\end{theorem}

\begin{proof}
\textbf{\ } Suppose that $y_{1},\ldots,y_{p}$ are elements in $\mathcal{M}$
that generate $\mathcal{M}$ as a von Neumann algebra. For every $0<\omega<1$,
there exists a family of noncommutative polynomials $\psi_{i}(x_{1}%
,\ldots,x_{n})$, $1\leq i\leq p$, such that
\[
\sum_{i=1}^{p}\Vert y_{i}-\psi_{i}(x_{1},\ldots,x_{n})\Vert_{2}^{2}<\left(
\frac{\omega}{4}\right)  ^{2}.
\]
For such a family of polynomials $\psi_{1},\ldots,\psi_{p}$, and
every $R>0$
there always exists a constant $D\geq1$, depending only on $R,\psi_{1}%
,\ldots,\psi_{n}$, such that
\[
\left (
\sum_{i=1}^{p}\Vert\psi_{i}(A_{1},\ldots,A_{n})-\psi_{i}(B_{1},\ldots
,B_{n})\Vert_{2}^{2}\right )^{1/2}\leq D
\Vert(A_{1},\ldots,A_{n})-(B_{1},\ldots ,B_{n})\Vert_{2},
\]
for all $(A_{1},\ldots,A_{n}),(B_{1},\ldots,B_{n})$ in $\mathcal{M}%
_{k}(\mathbb{C})^{ {n}}$, all $k\in\mathbb{N}$, satisfying $\Vert A_{j}\Vert,
\Vert B_{j}\Vert\leq R,$ for $1\leq j\leq n.$

For $R> 1, m$ sufficiently large, $\epsilon$ sufficiently small and $k$
sufficiently large, every $(H_{1},\ldots,H_{p},A_{1},\ldots,A_{n})$ in
$\Gamma_{R}(y_{1},\ldots,y_{p},x_{1},\ldots,x_{n};m,k,\epsilon)$ satisfies
\[
\left (\sum_{i=1}^{p}\Vert
H_{i}-\psi_{i}(A_{1},\ldots,A_{n})\Vert_{2}^{2}\right )^{1/2}\leq
\frac{\omega}{4} .
\]
It is obvious that such an $(A_{1},\ldots,A_{n})$ is also in $\Gamma_{R}%
(x_{1},\ldots,x_{n};m,k,\epsilon)$. On the other hand, by the
definition of the orbit covering number,  we know there exists a set
$\{ \mathcal{U}(B_{1}^{\lambda
},\ldots,B_{n}^{\lambda};\frac{\omega}{4D})
\}_{\lambda\in\Lambda_{k}}$ of $\frac{\omega}{4D}$-orbit-balls that
cover $\Gamma_{R}(x_{1},\ldots ,x_{n};m,k,\epsilon)$ with the
cardinality of $ \Lambda_{k} $ satisfying \
$|\Lambda_{k}|=\nu(\Gamma_{R}(x_{1},\ldots
,x_{n};m,k,\epsilon),\frac{\omega}{4D}).$ \ Thus \ for \ such \ $(A_{1}%
,\ldots,A_{n})$ in $\Gamma_{R}(x_{1},\ldots,x_{n};m,k,\epsilon)$, there exists
some $\lambda\in\Lambda_{k}$ and $W\in\mathcal{U}(k)$ such that
\[
\Vert(A_{1},\ldots,A_{n})-(WB_{1}^{\lambda}W^{\ast},\ldots,WB_{n}^{\lambda
}W^{\ast})\Vert_{2}\leq\frac{\omega}{4D}.
\]
It follows that
\[
\sum_{i=1}^{p}\Vert H_{i}-W\psi_{i}(B_{1}^{\lambda},\ldots,B_{n}^{\lambda
})W^{\ast}\Vert_{2}^{2}=\sum_{i=1}^{p}\Vert H_{i}-\psi_{i}(WB_{1}^{\lambda
}W^{\ast},\ldots,WB_{n}^{\lambda}W^{\ast})\Vert_{2}^{2}\leq\left(
\frac{\omega}{2}\right)  ^{2},
\]
for some $\lambda\in\Lambda_{k}$ and $W\in\mathcal U(k),$ i.e.,
\[
(H_{1},\ldots,H_{p})\in\mathcal{U}(\psi_{1}(B_{1}^{\lambda},\ldots
,B_{n}^{\lambda}),\ldots,\psi_{p}(B_{1}^{\lambda},\ldots,B_{n}^{\lambda
});\omega).
\]
Hence, by the definition of the free orbit-dimension, we get
\[
\begin{aligned}
0&\le \frak K(y_1,\ldots,y_p: x_1,\ldots,x_n;\omega,R)\le \inf_{m\in
\Bbb N, \epsilon>0}\limsup_{k\rightarrow \infty} \frac
{\log(|\Lambda_k|)}{-k^2\log\omega}\\
&= \inf_{m\in \Bbb N,
\epsilon>0}\limsup_{k\rightarrow \infty} \frac
{\log(\nu(\Gamma_R(x_1,\ldots, x_n ;m,k,\epsilon),\frac
\omega{4D}))}{-k^2\log\omega}\\
& = 0,
\end{aligned}
\]
since $\mathfrak{K}_{2}(x_{1},\ldots,x_{n})=0$. Therefore
$\mathfrak{K} (y_{1},\ldots,y_{p}:x_{1},\ldots,x_{n};\omega)=0$. Now
it follows from Lemma 2 that
\[
\mathfrak{K} (y_{1},\ldots,y_{p};\omega)=\mathfrak{K} (y_{1},\ldots
,y_{p}:x_{1},\ldots,x_{n};\omega)=0;
\]
whence $\mathfrak{K}_{2}(y_{1},\ldots,y_{p})=0$ and
$\mathfrak{K}_{2}( \mathcal{M})=0$.
\end{proof}

\begin{theorem}
If $\mathcal{M}$ is a hyperfinite von Neumann algebra with a tracial state
$\tau$, then $\mathfrak{K}_{2}( \mathcal{M})=0$.
\end{theorem}

\begin{proof}
When $\mathcal{M}$ is an abelian von Neumann algebra, the result follows from
\cite[Lemma 4.3]{V2}. Generally, it is a direct consequence of Proposition 1, that, for each $0<\omega<1,$%
\[
\nu\left(  \Gamma_{R}\left(  x_{1},\ldots,x_{n},m,\varepsilon,k\right)
,\omega\right)  =1
\]
whenever $m$ is sufficiently large and $\varepsilon$ is sufficiently small.
\end{proof}

The proof of next theorem, being a slight modification of that of
Theorem 1, will be omitted.

\begin{theorem}
Suppose that $\mathcal{M}$ is a finitely generated von Neumann algebra with a
tracial state $\tau$. Suppose that $\{\mathcal{N}_{i}\}_{i=1}^{\infty}$ is an
ascending sequence of von Neumann subalgebras of $\mathcal{M}$ such that
$\mathfrak{K}_{2}(\mathcal{N}_{i})=0$ for all $i\geq1$ and $\mathcal{M}%
=\overline{\cup_{i}\mathcal{N}_{i}}^{SOT}$. Then $\mathfrak{K}_{2}(
\mathcal{M})=0$.
\end{theorem}

\begin{definition}
A unitary matrix $U$ in $\mathcal{M}_{k}(\mathbb{C})$ is a \emph{Haar unitary
matrix} if $\tau_{k}(U^{m})=0$ for all $1\leq m<k$ and $\tau_{k}(U^{k})=1$.
\end{definition}

\noindent The proof of following lemma can be found in \cite{GS2} (
see also \cite{V4}).  For the sake of completeness, we also sketch
its proof here.

\begin{lemma}
Let $V_{1},V_{2}$ be two Haar unitary matrices in $\mathcal{M}_{k}%
(\mathbb{C})$. For every $\delta>0$, let
\[
\Omega(V_{1},V_{2};\delta)=\{U\in\mathcal{U}(k)\ |\ \Vert UV_{1}-V_{2}%
U\Vert_{2}\leq\delta\}.
\]
Then, for every $0< \delta<r$, there exists a set $\{Ball
(U_{\lambda}; \frac {4\delta} r)\}_{\lambda \in\Lambda}$ of
$\frac{4\delta}{r}$-balls in $\mathcal{U}(k)$ that cover
$\Omega(V_{1},V_{2};\delta)$ with the cardinality of $\Lambda$
satisfying $|\Lambda|\leq\left( \frac{3r}{2\delta}\right)
^{4rk^{2}}$.
\end{lemma}

\begin{proof}[Sketch of Proof]
\textbf{\ } Let $D$ be a diagonal unitary matrix, $diag(\lambda
_{1},\ldots
,\lambda _{k})$, where $\lambda _{j}$ is the $j$-th root of unity $1$. Since $%
V_{1},V_{2}$ are Haar unitary matrices, there exist $W_{1},W_{2}$ in $%
\mathcal{U}(k)$ such that $V_{1}=W_{1}DW_{1}^{\ast }$ and $%
V_{2}=W_{2}DW_{2}^{\ast }$. Let $\tilde{\Omega}(\delta )=\{U\in \mathcal{U}%
(k)\ |\ \Vert UD-DU\Vert _{2}\leq \delta \}.$ Clearly  $\Omega
(V_{1},V_{2};\delta )=\{W_{2}^{\ast }UW_{1}|U\in \tilde{\Omega}(\delta)\}$; whence $%
\tilde{\Omega}(\delta)$ and $\Omega (V_{1},V_{2};\delta )$ have the
same covering numbers.

Let $\{e_{st}\}_{s,t=1}^{k}$ be the canonical system of matrix units
of $ \mathcal{M}_{k}(\mathbb{C})$. Let
\[
\begin{aligned}
\mathcal  S_1 = span \{e_{st} \  | \  |\lambda_s-\lambda_t|< r \}
\qquad \mathcal S_2=M_k(\Bbb C) \ominus S_1.
\end{aligned}
\]%
For every $U=\sum_{s,t=1}^{k}x_{st}e_{st}$ in $\tilde{\Omega}(\delta
)$,
with $x_{st}\in \mathbb{C}$, let $T_{1}=\sum_{e_{st}\in \mathcal{S}%
_{1}}x_{st}e_{st}\in \mathcal{S}_{1}$ and $T_{1}=\sum_{e_{st}\in \mathcal{S}%
_{2}}x_{st}e_{st}\in \mathcal{S}_{2}$. But
\[\begin{aligned}
\delta ^{2}&\geq \Vert UD-DU\Vert _{2}^{2}=\sum_{s,t=1}^{k}|(\lambda
_{s}-\lambda _{t})x_{st}|^{2}\geq \sum_{e_{st}\in
\mathcal{S}_{2}}|(\lambda
_{s}-\lambda _{t})x_{st}|^{2}\\ &\geq r^{2}\sum_{e_{st}\in \mathcal{S}%
_{2}}|x_{st}|^{2}=r^{2}\Vert T_{2}\Vert _{2}^{2}.
\end{aligned}\]%
Hence $\Vert T_{2}\Vert _{2}\leq \frac{\delta }{r}$. Note that
$\Vert
T_{1}\Vert _{2}\leq \Vert U\Vert _{2}=1$ and $dim_{\mathbb{R}}{}\mathcal{S}%
_{1}\leq 4rk^{2}.$ By standard arguments on covering numbers, we
know
that $\tilde{\Omega}(\delta)$ can be covered by a set   $\{Ball(A^{\lambda };\frac{%
2\delta }{r})\}_{\lambda \in \Lambda }$ of $\frac{2\delta }{r}$-balls in $ \mathcal{M}_{k}(\mathbb{C})$ with $|\Lambda |\leq \left( \frac{3r}{2\delta }%
\right) ^{4rk^{2}}.$ Because $\tilde{\Omega}(\delta )\subset \mathcal{U}%
(k)$, after replacing   $A^{\lambda }$ by a unitary $U^{\lambda }$ in $%
Ball(A^{\lambda },\frac{2\delta }{r})$, we obtain that the set
$\{Ball (U_{\lambda };\frac{4\delta }{r} )\}_{\lambda \in \Lambda }$
of $\frac{4\delta }{r}$-balls in $\mathcal{U}(k)$ that cover
$\tilde{\Omega}(\delta )$ with the cardinality of $\Lambda $
satisfying $|\Lambda |\leq \left( \frac{3r}{2\delta }\right) ^{4r
k^{2}}$. The same result holds for $\Omega (V_{1},V_{2};\delta )$.
\end{proof}

\begin{definition}
Suppose that $\mathcal{M}$ is a diffuse von Neumann algebra with a tracial
state $\tau$. Then a unitary element $u$ in $\mathcal{M}$ is called a
\emph{Haar unitary} if $\tau(u^{m})=0$ when $m\neq0$.
\end{definition}

\begin{theorem}
Suppose $\mathcal{M}$ is a diffuse von Neumann algebra with a
tracial state $\tau$. Suppose $\mathcal{N}$ is a diffuse von Neumann
subalgebra of $\mathcal{M}$ and $u$ is a unitary element in
$\mathcal{M}$ such that $\mathfrak{K}_{2}(\mathcal{N})=0$ and
$\{\mathcal{N},u\}$ generates $\mathcal{M}$ as a von Neumann
algebra. If there exist Haar unitary elements $v_{1},v_{2},\ldots$
and   $w_{1},w_{2},\ldots$ in
$\mathcal{N}$ such that $\left\Vert v_{n}u-uw_{n} \right\Vert _{2}%
\rightarrow0$, then $\mathfrak{K}_{2}( \mathcal{M})=0$. In particular, if
there are Haar unitary elements $v,w$ in $\mathcal{N},$ such that $vu=uw$,
then $\mathfrak{K}_{2}\left(  \mathcal{M}\right)  =0.$
\end{theorem}

\begin{proof}
\textbf{\ }
Suppose that $\{x_{1},\ldots,x_{n}\}$ is a family of generators of
$\mathcal{N}$. Then we know that
$\{x_{1},\ldots,x_{n}, u\}$ is a family of generators of
$\mathcal{M}$.

For every $0<\omega<1$, $0<r<1$, there exist an integer $p>0$ and
two Haar unitary elements $v_{p}, w_{p}$ in $\mathcal{N}$ such that
\[
\|v_{p}u-uw_{p}\|_{2}< \frac{r\omega}{65}.
\]
Note that $\{x_{1},\ldots,x_{n}, v_{p},w_{p}\} $ is also a family of
generators of $\mathcal{N}$.

For $R>1 $, $m\in\mathbb{N}$, $\epsilon>0$ and $k\in\mathbb{N}$, by the
definition of the orbit covering number, there exists a set $\{\mathcal{U}%
(B_{1}^{\lambda},\ldots,B_{n}^{\lambda}, V^{\lambda},W^{\lambda};\frac
{r\omega}{64})\}_{\lambda\in\Lambda_{k}}$ of $\frac{r\omega}{64} $-orbit-balls
in $\mathcal{M}_{k}(\mathbb{C})^{n+2}$ that cover $\Gamma_{R}(x_{1}%
,\ldots,x_{n}, v_{p}, w_{p};m,k,\epsilon)$, where the cardinality of
$\Lambda$ satisfies
$|\Lambda_{k}|=\nu(\Gamma_{R}(x_{1},\ldots,x_{n}, v_{p},
w_{p};m,k,\epsilon),\frac{r\omega}{64}).$ When $m$ is sufficient
large, $\epsilon$ is sufficient small, by Proposition 1 we can
assume that
all $V^{\lambda},W^{\lambda}$ are Haar unitary matrices in $\mathcal{M}%
_{k}(\mathbb{C})$.

For $m$ sufficiently large and $\epsilon$ sufficiently small, when
$(A_{1},\ldots,A_{n}, V, W,U)$ is contained in
$\Gamma_{R}(x_{1},\ldots,x_{n}, v_{p}, w_{p},u;m,k,\epsilon)$ then,
by Proposition 1, there exists a unitary element $U_{1}$ in
$\mathcal{U}(k)$ so that
\[
\|U_{1}-U\|_{2} <\frac{r\omega}{64}\qquad\ \text{ and } \ \qquad\Vert
VU_{1}-U_{1}W\Vert_{2}<\frac{r\omega}{64}.
\]
It is easy to see that   $(A_{1},\ldots,A_{n},V, W)$ is also in $\Gamma_{R}(x_{1}%
,\ldots,x_{n}, v_{p},w_{p};m,k,\epsilon)$. Since $\Gamma_{R}(x_{1}%
,\ldots,x_{n}, v_{p}, w_{p};m,k,\epsilon)$ is covered by the set
$\{\mathcal{U}(B_{1}^{\lambda},\ldots,B_{n}^{\lambda},V^{\lambda},W^{\lambda
};\frac{r\omega}{64})\}_{\lambda\in\Lambda_{k}}$ of $\frac{r\omega}{64}
$-orbit-balls, there exist some $\lambda\in\Lambda_{k}$ and $X\in
\mathcal{U}(k)$ such that
\[
\Vert(A_{1},\ldots,A_{n},V,W)-(XB_{1}^{\lambda}X^{\ast},\ldots,XB_{n}%
^{\lambda}X^{\ast},XV^{\lambda}X^{\ast},XW^{\lambda}X^{\ast})\Vert_{2}%
\leq\frac{r\omega}{64}.
\]
Hence,
\[
\Vert
V^{\lambda}X^{\ast}U_{1}X-X^{\ast}U_{1}XW^{\lambda}\Vert_{2}=\Vert
XV^{\lambda}X^{\ast}U_{1}-U_{1}XW^{\lambda}X^{\ast}\Vert_{2}\leq\frac{r\omega
}{16}.
\]
Note that $V^{\lambda},W^{\lambda}$ were chosen to be Haar unitary
matrices in $\mathcal{M}_{k}(\mathbb{C})$. From Lemma 3, it follows
that there exists a set $\{Ball
(U_{\lambda,\sigma};\frac{\omega}{4})\}_{\sigma\in\Sigma_{k}}$ of
$\frac{\omega}{4}
$-balls in $\mathcal{U}(k)$ that cover $\Omega(V^{\lambda},W^{\lambda}%
;\frac{r\omega}{16})$ with $|\Sigma_{k}|\leq\left(
\frac{24}{\omega}\right) ^{4rk^{2}}$, i.e., there exists some
$U_{\lambda,\sigma}$ in $\{U_{\lambda
,\sigma}\}_{\sigma\in\Sigma_{k}}$ such that
\[
\Vert X^{\ast}U_{1}X-U_{\lambda,\sigma}\Vert_{2}=\Vert U_{1}-XU_{\lambda
,\sigma}X^{\ast}\Vert_{2}\leq\frac{\omega}{4}.
\]
Thus for such an $(A_{1},\ldots,A_{n},V, W,U)$ in
$\Gamma_{R}(x_{1},\ldots,x_{n},
v_{p}, w_{p},u;m,k,\epsilon)$, there exists some $(B_{1}^{\lambda}%
,\ldots,B_{n}^{\lambda},V^{\lambda},W^{\lambda})$ and $U_{\lambda,\sigma}$
such that
\[
\Vert(A_{1},\ldots,A_{n}, U)-(XB_{1}^{\lambda}X^{\ast},\ldots,XB_{n}^{\lambda
}X^{\ast}, XU_{\lambda,\sigma}X^{\ast})\Vert_{2}\leq\frac\omega2,
\]
for some $X\in\mathcal{U}(k)$, i.e.,
\[
(A_{1},\ldots,A_{n},
U)\in\mathcal{U}(B_{1}^{\lambda},\ldots,B_{n}^{\lambda},
U_{\lambda,\sigma};\omega).
\]
Hence, by the definition of the free orbit-dimension, we have shown
\[
\begin{aligned}
0\le  \frak K(x_1, \ldots,x_n, u: v_p,w_p;\omega,R)  &\le \inf_{m\in
\Bbb N, \epsilon>0}\limsup_{k\rightarrow \infty} \frac
{\log(|\Lambda_k||\Sigma_k|)}{-k^2\log\omega}\\
&\le \inf_{m\in \Bbb N, \epsilon>0}\limsup_{k\rightarrow \infty}
\left (\frac {\log(|\Lambda_k| )}{-k^2\log\omega} + \frac {\log
\left ( \frac {24} { \omega}\right )^{4rk^2}}{-k^2\log\omega}\right )\\
&\le 0+ 4r \cdot \frac {\log 24 -\log\omega}{-\log\omega},
\end{aligned}
\]
since $\mathfrak{K}_{2}(x_{1},\ldots,x_{n}, v_{p}, w_{p})\leq\mathfrak{K}%
_{2}(\mathcal{N})=0$. Thus, by Lemma 2,
\[
\begin{aligned}
0\le  \frak K(x_1, \ldots,x_n, u ;\omega)   =\frak K(x_1,
\ldots,x_n, u: v_p,w_p;\omega) \le  4r \cdot \frac {\log 24
-\log\omega}{-\log\omega}.
\end{aligned}
\]
Because $r$ is an arbitrarily small positive number, we have $\mathfrak{K}
(x_{1},\ldots,x_{n}, u;\omega)=0$; whence, $\mathfrak{K}_{2} (x_{1}%
,\ldots,x_{n}, u )=0$. By Theorem 1, $\mathfrak{K}_{2}(
\mathcal{M})=0$.
\end{proof}

\vspace{0.2cm} Using the results in \cite[Theorem 18]{DH}, the
preceding theorem can be easily extended as follows.

\begin{theorem}
Suppose $\mathcal{M}$ is a von Neumann algebra with a tracial state
$\tau$. Suppose $\mathcal{N}$ is a von Neumann subalgebra of
$\mathcal{M}$ and $a$ is an element in $\mathcal{M}$ such that
$\mathfrak{K}_{2}(\mathcal{N})=0$, and $\{\mathcal{N},a\}$ generates
$\mathcal{M}$ as a von Neumann algebra. If there exist two normal
operators $b_{1},b_{2}$ in $\mathcal{N}$ such that $b_{1}$,
$b_{2}$ have no common eigenvalues and $ab_{1}=b_{2}a$, then $\mathfrak{K}%
_{2}( \mathcal{M})=0$.
\end{theorem}

\begin{theorem}
Suppose $\mathcal{M}$ is a von Neumann algebra with a tracial state $\tau$.
Suppose $\mathcal{M}$ is generated by von Neumann subalgebras $\mathcal{N}%
_{1}$ and $\mathcal{N}_{2}$ of $\mathcal{M}$. If $\mathfrak{K}_{2}%
(\mathcal{N}_{1})=\mathfrak{K}_{2}(\mathcal{N}_{2})=0$ and $\mathcal{N}%
_{1}\cap\mathcal{N}_{2}$ is a diffuse von Neumann subalgebra of $\mathcal{M}$,
then $\mathfrak{K}_{2}( \mathcal{M})=0$.
\end{theorem}

\begin{proof}
\textbf{\ } Suppose that $\{x_{1},\ldots,x_{n}\}$ is a family of generators of
$\mathcal{N}_{1}$ and $\{y_{1},\ldots,y_{p}\}$ a family of generators of
$\mathcal{N}_{2}$. Since $\mathcal{N}_{1}\cap\mathcal{N}_{2}$ is a diffuse von
Neumann subalgebra, we can find a Haar unitary $u$ in $\mathcal{N}_{1}%
\cap\mathcal{N}_{2}$.

For every $R>1+\max_{1\leq i\leq n,1\leq j\leq p}\{\Vert x_{i}\Vert,\Vert
y_{j}\Vert\}$, $0<\omega<\frac{1}{2n}$, $0<r<1$ and $m\in\mathbb{N}$,
$\epsilon>0$, $k\in\mathbb{N}$, there exists a set $\{\mathcal{U}%
(B_{1}^{\lambda},\ldots,B_{n}^{\lambda},U_{\lambda};\frac{r\omega}%
{24R})\}_{\lambda\in\Lambda_{k}}$ of
$\frac{r\omega}{24R}$-orbit-balls in
$\mathcal{M}_{k}(\mathbb{C})^{n+1}$   covering
$\Gamma_{R}(x_{1},\ldots ,x_{n},u;m,k,\epsilon)$ with
$|\Lambda_{k}|=\nu(\Gamma_{R}(x_{1},\ldots
,x_{n},u;m,k,\epsilon),\frac{r\omega}{24R})$.

Also there exists a set
$\{\mathcal{U}(D_{1}^{\sigma},\ldots,D_{p}^{\sigma
},U_{\sigma};\frac{r\omega}{24R})\}_{\sigma\in\Sigma_{k}}$ of
$\frac{r\omega }{24R}$-orbit-balls in
$\mathcal{M}_{k}(\mathbb{C})^{p+1}$ that cover
$\Gamma_{R}(y_{1},\ldots,y_{p},u;m,k,\epsilon)$ with
$|\Sigma_{k}|=\nu (\Gamma_{R}(y_{1},\ldots,y_{p},u;
m,k,\epsilon),\frac{r\omega}{24R})$. When $m$ is sufficiently large
and $\epsilon$ is sufficiently small, by Proposition 1 we can assume
all $U_{\lambda}$, $U_{\sigma}$ to be Haar unitary matrices in
$\mathcal{M}_{k}(\mathbb{C})$.

For each $(A_{1},\ldots,A_{n},C_{1},\ldots,C_{p},U)$ in $\Gamma_{R}%
(x_{1},\ldots,x_{n},y_{1},\ldots,y_{p},u;m,k,\epsilon)$, we know
that $(A_{1},\ldots,A_{n},U)$ is contained in
$\Gamma_{R}(x_{1},\ldots,x_{n},u;m,k,\epsilon)$
and $(C_{1},\ldots,C_{p},U)$ is contained in $\Gamma_{R}(y_{1},\ldots,y_{p}%
,u;m,k,\epsilon)$. Note $\Gamma_{R}(x_{1},\ldots,x_{n},u;m,k,\epsilon)$ is
covered by the set $\{\mathcal{U}(B_{1}^{\lambda},\ldots,B_{n}^{\lambda
},U_{\lambda};\frac{r\omega}{24R})\}_{\lambda\in\Lambda_{k}}$ of
$\frac{r\omega}{24R}$-orbit-balls and $\Gamma_{R}(y_{1},\ldots,y_{p}%
,u;m,k,\epsilon)$ is covered by the set $\{\mathcal{U}(D_{1}^{\sigma}%
,\ldots,D_{p}^{\sigma},U_{\sigma};\frac{r\omega}{24R})\}_{\sigma\in\Sigma_{k}%
}$ of $\frac{r\omega}{24R}$-orbit-balls. Hence, there exist some $\lambda
\in\Lambda_{k}$, $\sigma\in\Sigma_{k}$ and $W_{1},W_{2}$ in $\mathcal{U}%
\left(  k\right)  $ such that
\[
\begin{aligned}
&\|(A_1,\ldots,A_n,U) - (W_1B_1^\lambda W_1^*, \ldots, W_1B_n^\lambda
W_1^*, W_1U_\lambda W_1^*)\|_2 \le \frac {r\omega}{24R}\\
& \|( C_1,\ldots,C_p, U)-(W_2D_1^\sigma W_2^*, \ldots, W_2D_p^\sigma
W_2^*, W_2U_\sigma W_2^*)\|_2\le \frac {r\omega}{24R}.
\end{aligned}
\]
Hence,
\[
\Vert
W_{2}^{\ast}W_{1}U_{\lambda}-U_{\sigma}W_{2}^{\ast}W_{1}\Vert_{2}=\Vert
W_{1}U_{\lambda}W_{1}^{\ast}-W_{2}U_{\sigma}W_{2}^{\ast}\Vert_{2}\leq
\frac{r\omega}{12R}.
\]
From our assumption that $U_{\lambda},U_{\sigma}$ are Haar unitary
matrices in $\mathcal{M}_{k}(\mathbb{C})$, by Lemma 3 we know that
there exists a set $\{Ball (U_{\lambda\sigma\gamma}; \frac{\omega}
{3R})\}_{\gamma\in\mathcal{I}_{k}}$ of $\frac{\omega}  {3R}$-balls
in $\mathcal{U}(k)$ that cover $\Omega(U_{\lambda},U_{\sigma
};\frac{r\omega}{12R})$ with the cardinality of $\mathcal{I}_{k}$
never exceeding $\left( \frac{18R}{\omega}\right)  ^{4rk^{2}}.$ Then
there exists some $\gamma\in\mathcal{I}_{k}$ such that $\Vert
W_{2}^{\ast}W_{1}-U_{\lambda\sigma\gamma}\Vert_{2}\leq\frac{\omega}{3R}$.
This in turn implies
\[
\begin{aligned}
\|(A_1,\ldots, A_n,C_1,\ldots,C_p, U) - &
(W_2U_{\lambda\sigma\gamma}B_1^\lambda
U_{\lambda\sigma\gamma}^*W_2^*, \ldots,
W_2U_{\lambda\sigma\gamma}B_n^\lambda
U_{\lambda\sigma\gamma}^*W_2^*, \\
& \quad\qquad \quad\quad   W_2D_1^\sigma W_2^*, \ldots,
W_2D_p^\sigma W_2^*,W_2U_\sigma W_2^* )\|_2\le  n\omega\end{aligned}
\]
for some
$\lambda\in\Lambda_{k},\sigma\in\Sigma_{k},\gamma\in\mathcal{I}_{k}$
and $W_{2}\in\mathcal{U}(k)$, i.e.,
\[
(A_{1},\ldots,A_{n},C_{1},\ldots,C_{p},U)\in\mathcal{U}(U_{\lambda\sigma
\gamma}B_{1}^{\lambda}U_{\lambda\sigma\gamma}^{\ast},\ldots,U_{\lambda
\sigma\gamma}B_{n}^{\lambda}U_{\lambda\sigma\gamma}^{\ast},D_{1}^{\sigma
},\ldots,D_{p}^{\sigma},U_{\sigma};2n\omega).
\]
Hence, by the definition of the free orbit-dimension we get
\[
\begin{aligned}
\frak K(x_1,\ldots,x_n,&y_1,\ldots,y_p,u;2n\omega, R)   \le \inf_{m\in
\Bbb N, \epsilon>0}\limsup_{k\rightarrow \infty} \frac
{\log(|\Lambda_k||\Sigma_k||\mathcal I_k|)}{-k^2\log(2n\omega)}\\
&\le \inf_{m\in \Bbb N, \epsilon>0}\limsup_{k\rightarrow
\infty}\left ( \frac {\log(|\Lambda_k|)} {-k^2\log(2n\omega)} +
\frac {\log(|\Sigma_k|)}{-k^2\log(2n\omega)} + \frac {\log(|\mathcal
I_k|)}{-k^2\log(2n\omega)}\right )\\
&\le 0+\inf_{m\in \Bbb N, \epsilon>0}\limsup_{k\rightarrow \infty}
\frac {\log \left  ( \frac {18R} \omega \right )^{4rk^2}
}{-k^2\log(2n\omega)}\\ &\le 4r\cdot \frac {\log (18R)-\log \omega
}{- \log(2n\omega)},
\end{aligned}
\]
since $\mathfrak{K}_{2}(N_{1})=\mathfrak{K}_{2}(N_{2})=0$. Since $r$
is an arbitrarily small positive number, we get that $\mathfrak{K}
(x_{1},\ldots,x_{n},y_{1},\ldots,y_{p},u;2n\omega, R)=0$; whence
$\mathfrak{K}_{2}(x_{1},\ldots,x_{n},y_{1},\ldots,y_{p},u)=0$. By
Theorem 1, $\mathfrak{K}_{2}( \mathcal{M})=0$.
\end{proof}

\section{Applications}

In this section, we   discuss a few applications of the results from
the last section. Let $L(F_{n})$ denote the free group factor on $n$
generators. By
Voiculescu's fundamental result in \cite{V2}, we know $\delta_{0}%
(L(F_{n}))\geq n$, where $\delta_{0}$ is Voiculescu's   free entropy
dimension. By combining Theorem 1, 2, 3, 4, 5 and 6,  we can easily
obtain the results in \cite{V3}, \cite{Ge2}, \cite{GS2}, \cite{GS1},
 and \cite{V4}. Here are a few sample improvements.

The following lemma can be proved using Theorem 5.3 of \cite {ChSm}.

\begin{lemma}
If $\mathcal{M}$ is a II$_{1}$ factor with property $\Gamma$  with
the tracial state $\tau$, then there are a hyperfinite II$_1$ factor
$\mathcal R$ and  a sequence $\left\{ u_{n}\right\}  $ of   Haar
unitary
elements of $\mathcal{R}$ such that%
\[
\left\Vert u_{n}x-xu_{n}\right\Vert _{2}\rightarrow0
\]
for every $x\in\mathcal{M}.$
\end{lemma}

\begin{corollary}
If $\mathcal{M}$ is a II$_{1}$ factor with property $\Gamma$, then
$\mathfrak{K}_{2}(\mathcal{M})=0.$
\end{corollary}

\begin{proof}
Choose a hyperfinite II$_1$ factor $\mathcal R$ and  a  sequence of
Haar unitary elements  $ u_{1},u_{2},\ldots $ in $\mathcal R$ such
that $\lim_{n\rightarrow \infty}\| xu_n-u_nx\|_2=4$ for every $x$ in
$\mathcal M$. Since $\mathcal{R}$ is hyperfinite,
$\mathfrak{K}_{2}\left( \mathcal{R}\right) =0.$ If $\left\{
v_{1},v_{2},\ldots\right\}  $ is a sequence of Haar unitaries that
generate $\mathcal{M},$ it inductively follows from Theorem 4
that, for each $n\geq1$%
\[
\mathfrak{K}_{2}\left(  \left(  \mathcal{R}\cup\left\{  v_{1},\ldots
,v_{n}\right\}  \right)  ^{\prime\prime}\right)  =0.
\]
Whence, by Theorem 3, $\mathfrak{K}_{2}\left(  \mathcal{M}\right)
=0.$
\end{proof}

A maximal abelian self-adjoint subalgebra (or, masa) $\mathcal{A}$ in a
II$_{1}$ factor $\mathcal{M}$ is called a \emph{Cartan subalgebra} if the
\emph{normalizer algebra} of $\mathcal{A},$
\[
\mathcal{N}_{1}\left(  \mathcal{A}\right)  =\left\{  u\in\mathcal{U}\left(
\mathcal{M}\right)  :u^{\ast}\mathcal{A}u\subset\mathcal{A}\right\}
^{\prime\prime}%
\]
equals $\mathcal{M}$.   We define $\mathcal{N}_{k+1}\left(
\mathcal{A}\right)  =\mathcal{N}_{1}\left(  \mathcal{N}_{k}\left(
\mathcal{A}\right)  \right)  $ for $k\geq1$, and
$\mathcal{N}_{\infty}\left( \mathcal{A}\right)  =\left(
\bigcup_{1\leq k<\infty}\mathcal{N}_{k}\left( \mathcal{A}\right)
\right)  ^{\prime\prime}.$ The following is a direct consequence of
Theorems 4 and 3.

\begin{corollary}
Suppose $\mathcal{M}$ is a type II$_{1}$ factor,
and $\mathcal{A}$ is a diffuse von Neumann subalgebra with $\mathfrak{K}%
_{2}\left(  \mathcal{A}\right)  =0$. If $\mathcal{M}=\mathcal{N}_{k}\left(
\mathcal{A}\right)  $ for some $k,1\leq k\leq\infty,$ then $\mathfrak{K}%
_{2}\left(  \mathcal{M}\right)  =0$, and $\delta_{0}\left(  \mathcal{M}%
\right)  \leq1$.
\end{corollary}

Many important applications of free entropy to finite von Neumann algebras
(nonprime factors, some II$_{1}$ factors with property $T$) are consequences
of a result of L. Ge and J. Shen \cite{GS2}, which states that if
$\mathcal{M}$ is a II$_{1}$ von Neumann algebra generated by a sequence of
Haar unitary elements $\{u_{i}\}_{i=1}^{\infty}$ in $\mathcal{M}$ such that
each $u_{i+1}u_{i}u_{i+1}^{\ast}$ is in the von Neumann subalgebra generated
by $\{u_{1},\ldots,u_{i}\}$ in $\mathcal{M}$, then $\delta_{0}(\mathcal{M}%
)\leq1$. This result is an easy consequence of Theorem 4. Here is a
sample of a result that is stronger.

\begin{corollary}
Suppose $\mathcal{M}$ is a factor of type II$_{1}$ that is generated by a
family $\left\{  u_{ij}:1\leq i,j<\infty\right\}  $ of Haar unitary elements
in $\mathcal{M}$ such that

\begin{enumerate}
\item for each $i,j$, $u_{i+1,j}u_{ij}u_{i+1,j}^{\ast}$ is in the von Neumann
subalgebra generated by $\{u_{1j},\ldots,u_{ij}\};$ and

\item for each $j\geq1$, $\left\{  u_{1j},u_{2j},\ldots\right\}
\bigcap\left\{  u_{1,j+1},u_{2,j+1},\ldots\right\}  \neq\varnothing.$
\end{enumerate}

Then $\mathfrak{K}_{2}(\mathcal{M})=0$, $\delta_{0}(\mathcal{M})\leq1$. Thus
$\mathcal{M}$ is not *-isomorphic to any $L(F(n))$ for $n\geq2$.
\end{corollary}

\begin{remark}
\textbf{\ } Many new examples can be obtained by using the preceding
corollary. For example, suppose that $G$ is a group generated by elements
$a,b,c$ such that $ab^{2}=b^{3}a$ and $ac^{2}=c^{3}a$. The group von Neumann
algebra associated with $G$ is a type II$_{1}$ factor, and the preceding
corollary easily implies that $\mathfrak{K}_{2}(L(G))=0$ and $\delta
_{0}(L(G))\leq1.$
\end{remark}

The next two corollaries follows directly from Corollary 3.

\begin{corollary}
Suppose $\mathcal M$ is a nonprime II$_1$ factor, i.e. $\mathcal
M\simeq \mathcal N_1\otimes \mathcal N_2$ for some II$_1$ subfactors
$ \mathcal N_1, \mathcal N_2$. Then
 $\mathfrak{K}_{2}( \mathcal{M})=0$, $\delta_{0}(
\mathcal{M})\leq1$. Thus $\mathcal{M}$ is not *-isomorphic to any
$L(F(n))$ for $n\geq2$.
\end{corollary}

\begin{corollary}
If $\mathcal{M}=L(SL(\mathbb{Z},2m+1))$ is the group von Neumann algebra
associated with $SL(\mathbb{Z},2m+1)$ (the special linear group with integer
entries) for $m\geq1$, then $\mathfrak{K}_{2}( \mathcal{M})=0$, $\delta_{0}(
\mathcal{M})\leq1$. Thus $\mathcal{M}$ is not *-isomorphic to any $L(F(n))$
for $n\geq2$.
\end{corollary}

\vspace{0.2cm}

In \cite{GePopa} L. Ge and S. Popa defined a type II$_{1}$ factor to
be \emph{weakly }$n$\emph{-thin}, if it contains hyperfinite
subalgebras $\mathcal{R}_{0},\mathcal{R}_{1}$ and $n$ vectors
$\xi_{1},\ldots,\xi_{n}$ in $L^{2}\left(  \mathcal{M},\tau\right)  $
such that $\mathcal{M}=\overline
{span}^{\left\Vert { \cdot }\right\Vert _{2}}\left(  \mathcal{R}_{0}\{\xi_{1}%
,\ldots,\xi_{n}\}\mathcal{R}_{1}\right)  .$ They showed that
$L(F_{m})$ is not weakly $n$-thin for $m> 2+2n$. Motivated by these
facts, we have the following definition.

\begin{definition}
A type II$_{1}$ factor $\mathcal{M}$ with the tracial state $\tau$ is
\emph{weakly }$\mathfrak{K}$\emph{-thin} (or, respectively, \emph{weakly }%
$n$\emph{-}$\mathfrak{K}$\emph{-thin}) if there exist von Neumann subalgebras
$\mathcal{N}_{0},\mathcal{N}_{1}$ of $\mathcal{M}$ with $\mathfrak{K}%
_{2}\left(  \mathcal{N}_{0}\right)  =\mathfrak{K}_{2}\left(  \mathcal{N}%
_{1}\right)  =0$ and a vector $\xi$ $($or, respectively, $n$ vectors $\xi
_{1},\ldots,\xi_{n}$$)$ in $L^{2}\left(  \mathcal{M},\tau\right)  $ such that
$\overline{span}^{\Vert\cdot\Vert_{2}}\left(  \mathcal{N}_{0}\xi
\mathcal{N}_{1}\right)  =L^{2}\left(  \mathcal{M},\tau\right)  $ $($or,
respectively, $\overline{span}^{\Vert\cdot\Vert_{2}}\mathcal{N}_{0}\{\xi
_{1},\ldots,\xi_{n}\}\mathcal{N}_{1}=L^{2}\left(  \mathcal{M},\tau\right)
$$)$.
\end{definition}

\begin{theorem}
Suppose that  $\mathcal{M}$ is a finitely generated \emph{weakly }$n$%
\emph{-}$\mathfrak{K}$\emph{-thin }type II$_{1}$ factor with a tracial state
$\tau$. Then $\mathfrak{K}_{1}( \mathcal{M})\leq1+2n$ and $\delta_{0}(
\mathcal{M})\leq2+2n.$ Thus $\mathcal{M}$ is not *-isomorphic to $L(F_{m})$
for $m> 2+2n$.
\end{theorem}

\begin{proof}
\textbf{\ } Suppose $x_{1},\ldots,x_{p}$ is a family of self-adjoint elements
in $\mathcal{M}$ that generate $\mathcal{M}$ as a von Neumann algebra. Note
there exist von Neumann subalgebras $\mathcal{N}_{0},\mathcal{N}_{1}$ of
$\mathcal{M}$ with $\mathfrak{K}_{2}(\mathcal{N}_{0})=\mathfrak{K}%
_{2}(\mathcal{N}_{1})=0$ and $n$ vectors $\xi_{1},\ldots,\xi_{n}$ in $L^{2}(
\mathcal{M},\tau)$ such that $\overline{span}^{\Vert\cdot\Vert_{2}}%
\mathcal{N}_{0}\{\xi_{1},\ldots,\xi_{n}\}\mathcal{N}_{1}=L^{2}( \mathcal{M}%
,\tau)$. We can choose self-adjoint elements $y_{1},y_{2},\ldots ,$
$y_{2n-1},y_{2n}$ in $\mathcal{M}$ to approximate
$Re\xi_{1},Im\xi_{1},\ldots, Re\xi_{n},Im\xi_{n}$, respectively.
Hence, for any positive $\omega<1$, there are a positive integer
$N$, elements $\{a_{i,j,l}\}_{1\leq i\leq p,1\leq j\leq N,1\leq
l\leq2n}$ in $\mathcal{N}_{0}$, $\{b_{i,j,l}\}_{1\leq i\leq p,1\leq
j\leq N,1\leq l\leq2n}$ in $\mathcal{N}_{1}$, and self-adjoint
elements $y_{1},\ldots,y_{2n}$ in $\mathcal{M}$ such that
\[
\sum_{i=1}^{p}\Vert x_{i}-\sum_{j=1}^{N}\sum_{l=1}^{2n}a_{i,j,l}y_{l}%
b_{i,j,l}\Vert_{2}^{2}\leq\left(  \frac{\omega}{8}\right)  ^{2}.
\]
Without loss of generality, we can assume that $\{a_{i,j,l}\}_{1\leq i\leq
p,1\leq j\leq N,1\leq l\leq n}$ generates $\mathcal{N}_{0}$ and $\{b_{i,j,l}%
\}_{1\leq i\leq p,1\leq j\leq N,1\leq l\leq n}$ generates $\mathcal{N}_{1}$ as
von Neumann algebras. Otherwise we should add generators of $\mathcal{N}_{0}$,
$\mathcal{N}_{1}$ into the families.

Let $a$ be $\max_{1\leq i\leq p}\{\Vert x_{i}\Vert_{2}\}+2$. From now on the
sequence $z_{1},\ldots,z_{s},\ldots,z_{t}$ is denoted by $(z_{s}%
)_{s=1,\ldots,t}$ or $(z_{s})_{s}$ if there is no confusion arising from the
range of index, where $z_{s}$ is an element in $\mathcal{M}$ or a matrix in
$\mathcal{M}_{k}(\mathbb{C})$.

For $R>a$, define mapping $\psi:( \mathcal{M}_{k}(\mathbb{C})^{N})^{2n}%
\times\mathcal{M}_{k}(\mathbb{C})^{2n}\times( \mathcal{M}_{k}(\mathbb{C}%
)^{N})^{2n}\rightarrow\mathcal{M}_{k}(\mathbb{C})$ as follows,
\[
\psi((D_{j,l})_{jl},(E_{l})_{l},(F_{j,l})_{jl})=\sum_{j=1}^{N}\sum_{l=1}%
^{2n}D_{j,l}E_{l}L_{j,l}.
\]
Let $( \mathcal{M}_{k}(\mathbb{C}))_{R}$ be the collection of all $A$ in
$\mathcal{M}_{k}(\mathbb{C})$ such that $\Vert A\Vert\leq R$. Then there
always exists a constant $D>1$, not depending on $k$, such that
\begin{align}
  \Vert &\left(  \psi(  (A_{1,j,l}^{(1)})_{jl},(Y_{l})_{l},(B_{1,j,l}%
^{(1)})_{jl} ) ,\ldots,\psi( (A_{p,j,l}^{(1)})_{jl},(Y_{l})_{l},(B_{p,j,l}%
^{(1)})_{jl} ) \right) \\
&   \qquad \qquad -\left(  \psi( (A_{1,j,l}^{(2)})_{jl},(Y_{l})_{l},(B_{1,j,l}%
^{(2)})_{jl} ) ,\ldots,\psi( (A_{p,j,l}^{(2)})_{jl},(Y_{l})_{l},(B_{p,j,l}%
^{(2)})_{jl} ) \right)  \Vert_{2}\nonumber\\
&  \quad\leq D\Vert\left(  (A_{i,j,l}^{(1)})_{ijl},(B_{i,j,l}^{(1)}%
)_{ijl}\right)  -\left(  (A_{i,j,l}^{(2)})_{ijl},(B_{i,j,l}^{(2)}%
)_{ijl}\right)  \Vert_{2},\nonumber
\end{align}
for all
\[
\left\{  A_{i,j,l}^{(1)},Y_{l},B_{i,j,l}^{(1)},A_{i,j,l}^{(2)},B_{i,j,l}%
^{(2)}\right\}  _{i,j,l}\subset(\mathcal{M}_{k}(\mathbb{C}))_{R}\qquad\forall
k\in N.
\]
For $m$ sufficiently large, $\epsilon$ sufficiently small and $k$ sufficiently
large, if
\[
\left(  X_{1},\ldots,X_{p},(A_{i,j,l})_{ijl},(Y_{l})_{l},(B_{i,j,l}%
)_{ijl}\right)  \in\Gamma_{R}(x_{1},\ldots,x_{p},(a_{i,j,l})_{ijl},(y_{l}%
)_{l},(b_{i,j,l})_{ijl};k,m,\epsilon),
\]
then
\begin{align}
\Vert(X_{1},\ldots,X_{p})-  &  (\psi((A_{1,j,l})_{jl},(Y_{l})_{l}%
,(B_{1,j,l})_{jl}),\ldots,\psi((A_{p,j,l})_{jl},(Y_{l})_{l},(B_{p,j,l}%
)_{jl}))\Vert_{2}\\
& \qquad  =(\ \sum_{i=1}^{p}\Vert X_{i}-\sum_{j=1}^{N}\sum_{l=1}^{2n}A_{i,j,l}%
Y_{l}B_{i,j,l}\Vert_{2}^{2}\ )^{1/2}\leq  \frac{\omega}{8}
,\nonumber
\end{align}
and
\[
\begin{aligned}
&( (A_{i,j,l})_{ijl}) \in \Gamma_R(  (a_{i,j,l})_{ijl} ;
k,m,\epsilon),\quad \text {and} \quad (  (B_{i,j,l})_{ijl})  \in
\Gamma_R( (b_{i,j,l})_{ijl} ; k,m,\epsilon).
\end{aligned}
\]
On the other hand, from the definition of the orbit covering number,
it follows there exists a set
$\{\mathcal{U}((A_{ijl}^{\lambda})_{ijl};\frac{\omega}{16D})\}_{\lambda
\in\Lambda_{k}}$, or
$\{\mathcal{U}((B_{ijl}^{\sigma})_{ijl};\frac{\omega
}{16D})\}_{\sigma\in\Sigma_{k}}$, of
$\frac{\omega}{16D}$-orbit-balls that
cover $\Gamma_{R}((a_{i,j,l})_{ijl};k,m,\epsilon)$, or $\Gamma_{R}%
((b_{i,j,l})_{ijl};k,m,\epsilon)$ respectively, with
\[
|\Lambda_{k}|=\nu(\Gamma_{R}((a_{i,j,l})_{ijl};k,m,\epsilon),\frac{\omega
}{16D}), \quad|\Sigma_{k}|=\nu(\Gamma_{R}((b_{i,j,l})_{ijl};k,m,\epsilon
),\frac{\omega}{16D}).
\]
Therefore for such sequence $((A_{i,j,l})_{ijl},(B_{i,j,l})_{ijl})$, there
exist some $\lambda\in\Lambda_{k}$, $\sigma\in\Sigma_{k}$ and $W_{1},W_{2}$ in
$\mathcal{U}\left(  k\right)  $ such that
\begin{align}
\Vert\left(  (A_{i,j,l})_{ijl},(B_{i,j,l})_{ijl}\right)  -((W_{1}%
A_{i,j,l}^{\lambda}W_{1}^{\ast})_{ijl},(W_{2}B_{i,j,l}^{\sigma}W_{2}^{\ast
})_{ijl})\Vert_{2}\leq\frac{\omega}{8D}.
\end{align}
Thus, from (4.1), (4.2) and (4.3), it follows that
\begin{align}
&  \Vert(X_{1},\ldots,X_{p}) -\left(  \psi((W_{1}A_{1,j,l}^{\lambda}W_{1}%
^{*})_{jl},(Y_{l})_{l},(W_{2}B_{1,j,l}^{\sigma}W_{2}^{*})_{jl}),\right.
\\
&  \qquad\qquad\qquad\qquad\qquad\qquad\ldots\left.  ,\psi((W_{1}%
A_{p,j,l}^{\lambda}W_{1}^{*})_{jl},(Y_{l})_{l},(W_{2}B_{p,j,l}^{\sigma}%
W_{2}^{*})_{jl})\right)  \Vert_{2} \nonumber\\
&  \qquad\qquad=\left ( \ \sum_{1\leq i\leq p}\Vert X_{i}-\sum_{j=1}^{N}\sum_{l=1}%
^{2n}W_{1}A_{i,j,l}^{\lambda}W_{1}^{\ast}Y_{l}W_{2}B_{i,j,l}^{\sigma}%
W_{2}^{\ast}\Vert_{2}^{2}\ \right )^{1/2}\leq  \frac{\omega}{4}
.\nonumber
\end{align}
Hence
\begin{align}
\left (\sum_{1\leq i\leq p}\Vert W_{1}^{\ast}X_{i}W_{1}-\sum_{j=1}^{N}\sum_{l=1}%
^{2n}\left(
A_{i,j,l}^{\lambda}W_{1}^{\ast}Y_{l}W_{2}B_{i,j,l}^{\sigma }\right)
W_{2}^{\ast}W_{1}\Vert_{2}^{2}\right )^{1/2}\leq \frac{\omega}{4}.
\end{align}

By a result of Szarek, there exists a $\frac{\omega}{4ap}$-net
$\{U_{\gamma }\}_{\gamma\in_{k}}$ in $\mathcal{U}(k)$ that cover
$\mathcal{U}(k)$ with respect to the uniform norm such that the
cardinality of $\mathcal{I}_{k}$ does not exceed
$(\frac{4apC}{\omega})^{k^{2}}$, where $C$ is a universal constant.
Thus  $\Vert
W_{2}^{\ast}W_{1}-U_{\gamma}\Vert\leq\frac{\omega}{4ap},$ for  some
$\gamma\in\mathcal{I}_{k}$. Because of (4.5), we know
\begin{align}
\Vert\sum_{j=1}^{N}\sum_{l=1}^{2n}A_{i,j,l}^{\lambda}W_{1}^{\ast}Y_{l}%
W_{2}B_{i,j,l}^{\sigma}\Vert_{2}\leq\Vert X_{i}\Vert_{2}+\omega<a.
\end{align}
From (4.5) and (4.6), we have
\begin{align}
\left (\sum_{1\leq i\leq p}\Vert W_{1}^{\ast}X_{i}W_{1}-\left(  \sum_{j=1}^{N}%
\sum_{l=1}^{2n}A_{i,j,l}^{\lambda}W_{1}^{\ast}Y_{l}W_{2}B_{i,j,l}^{\sigma
}\right)  U_{\gamma}\Vert_{2}^{2}\right )^{1/2}\leq  \frac{\omega}{2} %
\end{align}
Define a linear mapping $\Psi_{\lambda\sigma\gamma}: \mathcal{M}%
_{k}(\mathbb{C})^{2n}\mathbb{\rightarrow}
\mathcal{M}_{k}\mathbb{(C})^{p} $ as follows;
\[
\begin{aligned}
\Psi_{ \lambda\sigma\gamma}   (S_1, \ldots, S_{2n})  = \left  (
\frac 1 2\sum_{j=1}^N\sum_{l=1}^{2n} \left (A_{i,j,l}^{\lambda}  S_l
B_{i,j,l}^{\sigma}\right )U_\gamma  +\left (\left
(A_{i,j,l}^{\lambda} S_l   B_{i,j,l}^{\sigma}\right )U_\gamma
\right)^*\right )_{i=1,\ldots,p}.\end{aligned}
\]
Let $\mathfrak{F}_{\lambda\sigma\gamma}$ be the range of $\Psi_{\lambda
\sigma\gamma}$ in $\mathcal{M}_{k}(\mathbb{C})^{p} $. It is easy to see that
$\mathfrak{F}_{\lambda\sigma\gamma}$ is a real-linear subspace of
$\mathcal{M}_{k}(\mathbb{C})^{p} $ whose real dimension does not exceed
$2nk^{2}$. Therefore the bounded subset
\begin{align}
\{(H_{1},\ldots,H_{p})\in\mathfrak{F}_{\lambda\sigma\gamma}\ |\ \Vert
(H_{1},\ldots,H_{p})\Vert_{2}\leq ap\}
\end{align}
of $\mathcal{M}_{k}(\mathbb{C}){^{p}}$ can be covered by a set $\{(H_{1}%
^{\lambda\sigma\gamma,\rho},\ldots H_{p}^{\lambda\sigma\gamma,\rho}%
)\}_{\rho\in\mathcal{S}_{k}}$ of $\omega$-balls with the cardinality
of $\mathcal{S}_{k}$ satisfying
$|\mathcal{S}_{k}|\leq(\frac{3ap}{\omega })^{2nk^{2}}.$ But we know
from (4.6) that
\begin{align}
& \left \|\left(  \frac1 2\sum_{j=1}^{N}\sum_{l=1}^{2n} \left(
A_{i,j,l}^{\lambda }W_{1}^{*} Y_{l} W_{2}B_{i,j,l}^{\sigma}\right)
U_{\gamma}+ \left(  \left( A_{i,j,l}^{\lambda}W_{1}^{*} Y_{l}
W_{2}B_{i,j,l}^{\sigma}\right)  U_{\gamma
}\right)  ^{*} \right)  _{i=1,\ldots, p}\right \|_{2} \\
& =\left (\sum_{i=1}^{p}\|\frac1 2\sum_{j=1}^{N}\sum_{l=1}^{2n}
\left( A_{i,j,l}^{\lambda}W_{1}^{*} Y_{l}
W_{2}B_{i,j,l}^{\sigma}\right)  U_{\gamma}+ \left(  \left(
A_{i,j,l}^{\lambda}W_{1}^{*} Y_{l} W_{2}B_{i,j,l}^{\sigma
}\right)  U_{\gamma}\right)  ^{*}\|_{2}^{2}\right )^{1/2}\nonumber\\
& < ap ,\nonumber
\end{align}
and from (4.7) we have
\begin{align}
\|(  &  W_{1}^{*}X_{1}W_{1}, \ldots, W_{1}^{*}X_{p}W_{1})- \Psi_{\lambda
\sigma\gamma}(W_{1}^{*} Y_{1} W_{2},\ldots,W_{1}^{*} Y_{2n} W_{2}%
)\|_{2}\\
&  =\|(W_{1}^{*}X_{1}W_{1}, \ldots, W_{1}^{*}X_{p}W_{1})-\nonumber\\
&  \qquad\left(  \frac1 2\sum_{j=1}^{N}\sum_{l=1}^{2n} \left(  A_{i,j,l}%
^{\lambda}W_{1}^{*} Y_{l} W_{2}B_{i,j,l}^{\sigma}\right)  U_{\gamma}+ \left(
\left(  A_{i,j,l}^{\lambda}W_{1}^{*} Y_{l} W_{2}B_{i,j,l}^{\sigma}\right)
U_{\gamma}\right)  ^{*} \right)  _{i=1,\ldots, p} \|_{2}\nonumber\\
&  \le\omega.\nonumber
\end{align}
Thus, from (4.8), (4.9) and (4.10), there exists some
$\rho\in\mathcal{S}_{k}$ such that
\[
\Vert(W_{1}^{\ast}X_{1}W_{1},\ldots,W_{1}^{\ast}X_{p}W_{1})-(H_{1}%
^{\lambda\sigma\gamma,\rho},\ldots
H_{p}^{\lambda\sigma\gamma,\rho})\Vert _{2} \leq2\omega.
\]
By the definition of the free orbit-dimension, we know that
\[
\begin{aligned}
\frak K&(x_1,\ldots, x_p:  (a_{ijl})_{ijl}, (y_l)_{l},
(b_{ijl})_{ijl};4\omega,R) \le \inf_{m\in \Bbb N,
\epsilon>0}\limsup_{k\rightarrow \infty} \frac
{\log(|\Lambda_k||\Sigma_k||\mathcal I_k||\mathcal
S_k|)}{-k^2\log(4\omega)}\\
&\le \inf_{m\in \Bbb N, \epsilon>0}\limsup_{k\rightarrow \infty}
\left (\frac {\log |\Lambda_k|  }{-k^2\log(4\omega)}+\frac {\log
|\Sigma_k|  }{-k^2\log(4\omega)}+\frac {\log(\frac {4apC}
\omega)^{k^2}(\frac {3ap} \omega)^{2nk^2} )}{-k^2\log(4\omega)}
\right
)\\
&= 0+0 + \frac {\log( 4\cdot (3ap)^{2n}\cdot apC)-(2n+1)\log \omega
} {- \log(4\omega)},
\end{aligned}
\]
since
$\mathfrak{K}_{2}(\mathcal{N}_{0})=\mathfrak{K}_{2}(\mathcal{N}_{1})=0$.
Thus, by Lemma 2
\[
\begin{aligned}
0\le \frak K (x_1,\ldots, x_p ;4\omega )&=\frak K (x_1,\ldots, x_p:
(a_{ijl})_{ijl}, (y_l)_{l}, (b_{ijl})_{ijl};4\omega) \\ &\le  \frac
{\log( 4\cdot (3ap)^{2n}\cdot apC)-(2n+1)\log \omega } {-
\log(4\omega)}.
\end{aligned}
\]
By the definition of the free orbit-dimension, we obtain
\[
\mathfrak{K}_{1}(x_{1},\ldots,x_{p}) \leq\limsup_{\omega\rightarrow0}
\frac{\log( 4\cdot(3ap)^{2n}\cdot apC)-(2n+1)\log\omega} {- \log(4\omega)}
\leq1+2n.
\]
Hence, $\mathfrak{K}_{1}( \mathcal{M})\leq1+2n$ \ and \ $\delta_{0}(
\mathcal{M})\leq2+2n$.
\end{proof}

\vspace{0.2cm}

\begin{remark}
The mapping $a\mapsto a^{\ast}$ extends from $\mathcal{M}$ to a unitary map on
$L^{2}\left(  \mathcal{M},\tau\right)  ,$ so for $\xi\in L^{2}\left(
\mathcal{M},\tau\right)  ,$ it makes sense to talk about ${Re}\xi=\left(
\xi+\xi^{\ast}\right)  /2$ and ${Im}\xi=\left(  \xi-\xi^{\ast}\right)  /2i.$
In particular, it makes sense to talk about self-adjoint elements of
$L^{2}\left(  \mathcal{M},\tau\right)  .$ If we have $\overline{span}%
^{\Vert\cdot\Vert_{2}}\mathcal{N}_{0}\{\xi_{1},\ldots,\xi_{n}\}\mathcal{N}%
_{1}=L^{2}\left(  \mathcal{M},\tau\right)  $ with
$\xi_{1},\ldots,\xi_{n}$ self-adjoint elements in $L^{2}\left(
\mathcal{M},\tau\right)  ,$ the proof of Theorem 7 yields
$\mathfrak{K}_{1}(\mathcal{M})\leq1+n$ and $\delta_{0}(
\mathcal{M})\leq2+n.$
\end{remark}

Combining Theorem 7 and the preceding remark with Theorem 3, we have
the following corollaries (see also   \cite{Ge1} and \cite{GePopa}).

\begin{corollary}
$L(F_{n})$ has no simple maximal abelian self-adjoint subalgebra for $n\geq4 $.
\end{corollary}

\begin{corollary}
$L(F_{n})$ is not a $\mathfrak{K}$-thin factor for $n\geq4$.
\end{corollary}

\begin{remark}
\textbf{\ } Another corollary of Theorem 7 is as follows. Suppose
$\mathcal{M}$ is a II$_{1}$ factor with a tracial state $\tau$.
Suppose that $\mathcal{N}$ is a subfactor of $\mathcal{M}$ with
finite index, i.e., $[ \mathcal{M}:\mathcal{N}]=r<\infty$. If
$\mathfrak{K}_{2}(\mathcal{N})=0$, \ then $\mathfrak{K}_{1}(
\mathcal{M})\leq2[r]+3$ and $\delta_{0}( \mathcal{M})\leq2[r]+4$
where $[r]$ is the integer part of $r$.
\end{remark}

\end{document}